\documentclass{amsart}

\usepackage{amssymb,amsmath}

\usepackage[T1]{fontenc}
\usepackage[latin1]{inputenc}
\usepackage[dvips]{graphicx}
\usepackage{enumerate}

 \newtheorem{lemma}{Lemma}
 \newtheorem{theorem}{Theorem}
 \newtheorem{definition}{Definition}
 \newtheorem{example}{Example}
 
 \newtheorem{remark}{Remark}

\newcommand{\nto}{n\to\infty}

\newcommand{\eqd}{\stackrel{d}{=}}

\newcommand{\bbr}{\mathbb{R}}
\newcommand{\bbn}{\mathbb{N}}

\newcommand{\lln  }{{\rm lln }}

\begin{document}
\title{On peeling procedure applied to a Poisson point process}
\author{Yu. Davydov${}^{1}$, A.V. Nagaev ${}^{2}$ and A. Philippe${}^{3}$ 
}
\address{ ${}^{1}$  Laboratoire Paul Painlev\'e  \\
  Universit\'e de Lille 1 \\  
Batiment  M2, 59655 Villeneuve d'Ascq Cedex, France \\
  \vskip 0.1cm
  ${}^{2}$Faculty of Mathematics and Computer Science, \\ Nicolaus Copernicus University, Chopina,  Torun, Poland\\
  \vskip 0.1cm
${}^{3}$ Laboratoire de math\'ematiques Jean Leray \\ 
Universit\'e de Nantes \\
  2 rue de la Houssini\`ere,  44322 Nantes Cedex 3, France
 } 

\maketitle

\begin{flushright}
  \begin{quote}
Main results of this paper were obtained together with Alexander Nagaev,
with whom the first author had collaborated for more than 35 years, until
Alexander's tragic death in 2005. Since then, we have gathered strength
and finalised this paper, strongly feeling Alexander's absence - our
memories of him will stay with us forever. 
\end{quote}
\end{flushright}

\begin{abstract}
    In the focus of our attention is the asymptotic properties of the
  sequence of convex hulls which arise as a result of a peeling
  procedure applied to the convex hull generated by a Poisson point
  process. Processes of the considered type are tightly connected with
  empirical point processes and stable random vectors. Results are given about the limit
  shape of the convex hulls in the case of a discrete spectral measure. We give some numerical experiments to
  illustrate the peeling procedure for a larger class of Poisson point
  processes.\\   
\textbf{ Keywords :}  Control measure ; convex hull  ;  limiting shape ; 
  peeling ; Poisson point processes ;  stable vectors.
\end{abstract}

\section{Introduction}

Consider a Poisson point process (p.p.p.) $\pi=\pi_{\alpha,\nu}$ with
points scattered over $\bbr^d.$  Identify
$\bbr^d\setminus\{{\mathbf 0}\}$ with $\bbr_+^1\times S^{d-1}$
where $S^{d-1}$ is the unit sphere. Assume that the intensity measure of this process 
$\mu$ is of the
form \begin{equation}\label{CONTROL} \mu=\theta\times\nu \end{equation} where 
\begin{enumerate}\item 
  $\theta$ is the absolutely continuous measure on $\bbr_+^1$
  determined by the density function
$$\frac{{\rm d}\theta}{{\rm d}\lambda}(r)=\alpha r^{-\alpha-1},\ r>0,$$
$\lambda$ is the Lebesgue measure in $\bbr^1,\ \alpha>0$ is a parameter while
\item $\nu$, called the spectral density,  is a bounded measure on the $\sigma-$algebra $\mathcal{
  B}_{S^{d-1}}$ of the Borel subsets of $S^{d-1}.$ Without loss of generality we assume that
$\nu(S^{d-1})=1.$ We denote  by $S_{\nu}$ the support of $\nu$. 
\end{enumerate}

The representation (\ref{CONTROL}) means that for any Borel
$A\subset\bbr^1$ and $\ E\subset S^{d-1}$
$$\mu\left\{x\; \Big|\ |x|\in A,\ e_x\in E\right\}=\theta(A)\nu(E).$$
where $e_x=|x|^{-1}x$, for all $x\in \bbr^d$-.

We assume that the Poisson point process $\pi_{\alpha,\nu}$  is non-unilateral. It means that
$\nu$ is supported by a set $S_\nu\subset S^{d-1}$ such that the cone 
$\mathrm{cone}(S_\nu)$ generated by $S_\nu$ coincides with $\bbr^d.$

Let $B({\mathbf 0},r)$ denote the ball of a radius $r$ centred at
the origin. It is easily seen that for any $\delta>0$
$$\mu\left(B({\mathbf 0},\delta)\right)=\infty,\qquad
\mu\left(\bbr^d\setminus B({\mathbf
0},\delta)\right)<\infty.$$ It implies that
with probability 1 in any neighbourhood of the origin there are
infinitely many points of $\pi$ while $\pi(\bbr^d\setminus B({\mathbf
0},\delta) )$ is finite.

The interest to the point processes controlled by (\ref{CONTROL})
is explained by the following facts.

Let $\xi^{(1)},\ \xi^{(1)},\ \dots,\ \xi^{(n)}$ be independent
copies of a random vector $\xi$ such that the function ${\rm
P}\{|\xi|>r\}$ regularly varies as $r\to\infty$ with the exponent
$-\alpha$ and the measures $(\nu_r)_{r}$ defined by 
$$\nu_r(E)={\rm P}\{e_{\xi}\in E \;  : \;  |\xi|>r\},\ E\in \mathcal{
B}_{S^{d-1}},$$ weakly converge  to $\nu$ on $B^c({\mathbf
0},\tau) $ for any $\tau>0$.

Consider the empirical point process $\beta_n$ generated by
$\xi^{(1)},\ \xi^{(1)},\ \dots,\ \xi^{(n)}$ or, more precisely, by
the random set
$$
\varsigma_1^{(n)}=\{b_n^{-1}\xi^{(1)},\dots,b_n^{-1}\xi^{(n)}\}=\{\hat{\xi}^{(1)},\dots,\hat{\xi}^{(n)}\}$$
where
$$b_n=\inf{}\left\{r \; : \;   n{\rm P}\left( e_{\xi}\in E \, \big|\, \ |\xi|>r \right) \le
  1\right\}.$$ It is easily seen that the point process $\beta_n$
weakly converges to $\pi=\pi_{\alpha,\nu}$ (see e.g. \cite{Resnik}, Prop.
3.21). Thus, each $\pi_{\alpha,\nu}$ is a weak limit of a sequence of the
empirical processes.

It can be easily established that $\pi_{\alpha,\nu}$ admits the following
representation 
\begin{equation}
  \label{eq:pp-rep-expo}
  \pi_{\alpha,\nu} =\sum_{k=1}^{\infty}
  \delta_{\{\Gamma_k^{-1/\alpha} \epsilon_k\}}
\end{equation}
where 
\begin{itemize}
\item $\Gamma_k= \displaystyle\sum_{i=1}^{k} \gamma_k$, the sequence $(\gamma_k)_{k\in\mathbb{N}}$ is a sequence
  of i.i.d. random variable with common exponential distribution with
  mean equal to one.
\item  $(\epsilon_k)_{k\in\mathbb{N}}$ is a sequence
  of i.i.d. with common distribution $\nu$
\item  $(\gamma_k)_{k\in\mathbb{N}}$ and  $(\epsilon_k)_{k\in\mathbb{N}}$ are supposed to be independent
  \end{itemize}
(see e.g. \cite{davetal}).
It is worth recalling that the point processes $\pi_{\alpha,\nu}$
naturally arise within the framework of the theory of stable
distributions. For example, if $\alpha\in(0,1)$ then the series
$\zeta=\sum_{x^{(j)}\in \pi_{\alpha,\nu} } x^{(j)}$ converges a.s.. Furthermore,
$\zeta$ has the $d-$dimensional stable distribution
( see e.g. \cite{lepage}, \cite{SamoTaq}, \cite{davetal}).

In the focus of our interest is the sequence of convex hulls that
arise from the peeling procedure introduced in \cite{yusacha}.  In what
follows by $C(A),\ A\subset\bbr^d,$ we denote the convex hull
generated by $A.$ Let $C$ be a convex set. By ${\rm ext}\ C$ we
denote the set of the extreme points of $C.$ If $C$ is a convex
polyhedron then ${\rm ext\ C}$ is the finite set of its vertices.

It convenient to start with the binomial process $(\varsigma_1^{(n)})$. Let $C_1^{(n)}=C(\varsigma_1^{(n)}).$ If the measure ${\rm P}_\xi$, the
distribution of $\xi$,  has no atoms then a.s. $C_1^{(n)}$ is a polyhedron and,
furthermore, $$\varsigma_1^{(n)}\cap\partial\ C_1^{(n)}={\rm ext}\
C_1^{(n)}.$$ Define $$\varsigma_2^{(n)}=\varsigma_1^{(n)}\setminus{\rm ext}\
C_1^{(n)},\ C_2^{(n)}=C(\varsigma_2^{(n)}),$$ then
$$\varsigma_3^{(n)}=\varsigma_2^{(n)}\setminus{\rm ext}\ C_2^{(n)},\
C_3^{(n)}=C(\varsigma_3^{(n)})$$ and so on. Obviously, the sequence of
the so-built non empty convex hulls $C_1^{(n)},\
C_2^{(n)},\dots,C_k^{(n)},\dots$ is finite and its length is
random.

\begin{definition}
    \label{LTRL} We say that the underlying distribution ${\rm P}_\xi$
  and corresponding to it spectral measure $\nu$ are {\em
    non-unilateral\/} if the minimal closed cone containing $S_{\nu}$
  coincides with $\bbr^d.$
\end{definition}

If the underlying distribution  ${\rm P}_\xi$ is
non-unilateral then for any fixed $k$ 
$${\mathbf 0}\in {\rm int}\ C_k^{(n)},\qquad
\inf\left\{|x| \Big|\ x\in\partial C_k^{(n)}\right\}>0\qquad
\mbox{a.s.} $$ for $n$  sufficiently large.  In  \cite{yusacha}, it was shown that if $0<\alpha<2$ and
$\nu$ is non-unilateral then $(C_0^{(n)},\ C_1^{(n)},\ \dots,\
C_k^{(n)}),$ as $\nto,$ converge in distribution to  $(C_0,\
C_1,\ \dots,\ C_k)$ for any fixed $k.$  Consequenlty, the sequence $\#\{\mathrm{ext} C_k^{(n)}\}$ is bounded in probability as $n$ tends to infinity . 

In order to learn how
$C_k^{(n)}$ relates to $C_k$ when $k=k_n\to\infty$ we need, first,
to learn how $C_k$ behaves as $k\to\infty.$ It should be noted
that $C_k^{(n)}$ can be regarded as the multi-dimensional analogue
of the
 order statistics. So, the asymptotic properties
of $C_{k_n}^{(n)}$ are of great interest from the view-point of
mathematical statistics. 

\vskip .4cm 

We generalize now the construction of the peeling sequence to infinite set. \\
Let $\varsigma=\varsigma_1$ denote the set of points of $\pi$ or, in other
words, let $\varsigma$ support the random measure $\pi.$ We may apply to
$\varsigma$ the same peeling procedure as in the case of the finite set
$\varsigma_1^{(n)}.$ As a result we obtain the sequence of sets $\varsigma_1,\
\varsigma_2,\dots,$ the sequence of their convex hulls $C_1,\ C_2,\dots$
and the sets of extreme points ${\rm ext}\ C_k,\ k=1,2,\dots\ .$
Furthermore, a.s. $\varsigma_{k+1}=\varsigma_k\setminus{\rm ext}\ C_k.$ If
$\nu$ is non-unilateral then ${\mathbf 0}\in C_k$ a.s. for all
$k.$ Furthermore, $C_k$ is a.s. a polyhedron and $$\varsigma_k\cap\partial\ C_k={\rm ext}\ C_k.$$

Intuitively, we expect that the asymptotic behaviour of $C_n$ is
rather regular. 
It is convenient to state our basic conjecture in the following
way: 
\begin{quote}
\textbf{Denote}
 \begin{equation}\label{XXXX}\hat{C}_n=\rho_n^{-1}C_n, \  \text{where} \  \rho_n=\max_{x\in C_n}\ |x|.\end{equation} 
\textbf{
  If $\nu$ is non-unilateral then
there exists a non-random set $\hat{C}$ such that
$$\lim_{\nto}d_H(\hat{C}_n,\hat{C})=0\qquad a.s.$$}
\end{quote}
\begin{definition}
  $\hat{C}$ (if it exists) is called the {\em limit shape\/} of the
sequence $\hat{C}_n$.
\end{definition}

It is easy to show that if such $\hat{C}$ exists then it is
certainly non-random. Indeed, labelling the points of $\varsigma$ in the
descending order of their distances from the origin we obtain a
sequence $x^{(1)},x^{(2)},x^{(3)},\dots$ such that a.s.
$$|x^{(1)}|>|x^{(2)}|>|x^{(3)}|>\cdots.$$
It is worth noting that the joint distribution of
$|x^{(1)}|,|x^{(2)}|,\dots,|x^{(n)}|$ is absolutely continuous
with the density of the form \begin{equation}\label{DENSITY}
p_n(r_1,r_2,\dots,r_n)=\alpha^n(r_1r_2\cdots r_n)^{-\alpha-1}
e^{-\nu(S^{d-1})r_1^{-\alpha}}\mathbb{I}_{\{ r_1>r_2>\cdots>r_n\}}. \end{equation}

Let $(\eta,\ \varepsilon),\ (\eta_1,\ \varepsilon^{(1)}),\ (\eta_2,\ \varepsilon^{(2)}),\dots$
be i.i.d. with common distribution 
$${\rm P}\{\eta>r,\ \varepsilon\in E\}=e^{-r}\nu(E).$$

According to (\ref{eq:pp-rep-expo}), we have  

\begin{equation}
\{x^{(j)}\}_{j=1}^{\infty}\eqd\{\varepsilon^{(j)}(\eta_1+\cdots+\eta_j)^{1/\alpha}\}_{j=1}^{\infty}, 
\label{ppsim}
\end{equation}
which implies that the event $\{\lim_{\nto}\hat{C}_n\
\mbox{exists}\}$ belongs to the $\sigma$-algebra $\mathcal{I}$ of the
events invariant with respect to all finite permutations of the
random vectors  $(\eta_1,\ \varepsilon^{(1)}),\ (\eta_2,\
\varepsilon^{(2)}),\dots.$ By the Hewitt-Savage zero-one law $\mathcal{I}$ is
trivial. Since the limit set $\hat{C}=\lim_{\nto}\hat{C}_n$ is
$\mathcal{I}$-measurable we conclude that  $\hat{C}$ is constant with
probability 1.

Now we give a first example where the existence of the limit shape is proved. 
\begin{example}
  \label{EXMPL} Let $S_{\nu}$
consist of $d+1$ unit vectors $e^{(1)},\dots,e^{(d+1)}$ such that  
${\rm cone}( e^{(1)},\dots,e^{(d+1)})$ coincides with $\bbr^d.$ 
$\pi_{{\alpha,\nu }}$ is decomposed on $d+1$ one-dimensional independent p.p.p. of the form  $\left( x_{k}^{(i)}=|x_{k}^{(i)}|e^{(i)} \right)$. 
Since $\nu$ is non-unilateral the points $x_{k}^{(i)},\
i=1,2,\dots,d+1,$ serve as vertices of $C_k,\ k=1,2,\dots\ .$
Moreover,  $|x_{k}^{(i)}|\eqd(\nu_i)^{1/\alpha}(\eta_1+\cdots+\eta_k)^{1/\alpha}$ and $\rho_nn^{1/\alpha}\to t^+=\max_{1\le i\le d+1}\ (\nu_i)^{1/\alpha},\ a.s. $. Then  the limit shape 
 $\hat{C}$ is the convex polyhedron with vertices
$v^{(i)}=(t_i/t^+)e^{(i)}$  and $t_i=(\nu_i)^{1/\alpha},\
i=1,2,\dots,d+1.$  
\end{example}

If $\#(S_\nu)>d+1$ then the situation becomes much more complicated.
Theorem \ref{SHAPE} and \ref{th2} proved below deal with a case where a non-unilateral $\nu$ is supported  by a finite number of unit vectors.

Intuitively, we expect that, say, in case of $\nu$ uniformly
distributed over $S^{d-1}$ the unit ball arises as the limit
shape. However, it is not easy at all to prove this formally. The
authors tried to verify the credibility of this conjecture using
the Monte Carlo simulation. Obviously, the representation
(\ref{eq:pp-rep-expo})  provides a basis for such a simulation. The results of
simulation presented below make this conjecture very credible.

It should be emphasised that the basic goal of the present paper
is to draw attention to new and interesting problems of stochastic 
geometry. So far, little or nothing is
known about the peels no matter what point process they concern.

The paper is organised as follows. In Section 2, we obtained a partial
result on the limit shape of the convex hulls $C_k(\pi_{\alpha,\nu})$ when 
the spectral measure of the process $\pi_\alpha$ is atomic.   
Section 4 contains some numerical experiments.

\section{Almost sure convergence of the peeling}

In this section we assume that the spectral measure $\nu$ of the
process $\pi_{\alpha,\nu}$ is atomic, i.e. it is supported by a finite
number of the points $e^{(1)},\ldots,e^{(l)}$ belonging to the
unit sphere $S^{d-1}$. 

Furthermore, it is also assumed that ${\rm
cone}\{e^{(1)},\ldots,e^{(l)}\} = \bbr^d.$ Denote by
$\nu_i=\nu(\{e^{(i)}\}),\ i=1,2,\dots,l,$ the atoms of $\nu.$

It implies that the considered  point process is a superposition
of the one-dimensional independent Poisson point processes defined on the rays
$$\mathcal{ L}_i=\{x|\ x=te^{(i)},\ t>0\}, \; i = 1, \ldots,l.$$ If a
Borel set $A\subset\mathcal{ L}_i$ then
\begin{equation}\label{CNTRLMSR}\mu(A)=\nu_i\alpha\int_{A}r^{-\alpha-1}{\rm d}r.\end{equation}

\begin{definition}
Let  $A=\{a^{(1)},\dots,a^{(m)}\}$ be a finite set where 
$a^{(i)}\in\bbr^d$  for $i=1,\dots,m,$ and $m\ge d+1.$  The set
$A$ is {\em extreme} if $${\rm ext}\ C(A)=A.$$
\end{definition}


\begin{theorem}
  \label{SHAPE} Let  $C_k(\pi_{\alpha,\nu})$ be the $k$-th convex
  hull of the Poisson point process $\pi_{\alpha,\nu}$. Denote by $C_{\infty}$ the convex hull generated
  by $A=\{\nu_1^{1/\alpha}e^{(1)},\ldots,\nu_l^{1/\alpha}e^{(l)}\}.$ 

If $A$
  is extreme then as $k\to\infty$ 
  \begin{equation}
\label{RATE1}
  d_H\left(k^{1/\alpha}C_k(\pi_{\alpha,\nu}),\; C_{\infty}\right) =
  O\left(\sqrt{\frac{\lln{k}}{k}}\right)\ a.s.
\end{equation}
  where $\lln k
  = \ln\ln k.$ 
 The polyhedron $C_{\infty}$ determines the limit shape
 of the convex hulls $C_k(\pi_{\alpha,\nu})$.
\end{theorem}

\begin{remark} If $\sigma$ is uniformly distributed over its support in the
sense that $\nu_i = l^{-1}$, then the total number of the vertices
of $C_{\infty}$ equals $l.$ Furthermore, they lie on the sphere of
the radius  $l^{-1/\alpha}.$ Loosely speaking, the convex hulls
$C_k(\pi_{\alpha,\nu})$ are getting round as $k\to\infty.$ \end{remark}

If the condition $A$ is extreme is omitted, we can state the following
result~:  
\begin{theorem}\label{th2}  Let $C_k(\pi_{\alpha,\nu})$ be the
$k$-th convex hull of $\pi_{\alpha,\nu}$. Denote by $C_{\infty}$ the convex
hull generated by
$A=\{\nu_1^{1/\alpha}e^{(1)},\ldots,\nu_l^{1/\alpha}e^{(l)}\}.$  Then as $k\to\infty$ 
 \begin{equation}\label{CV22}
d_H\left(k^{1/\alpha}C_k(\pi_{\alpha,\nu}),\; C_{\infty}\right) \to \; 0 
\; a.s. \end{equation}
\end{theorem} 
\begin{remark} \label{az}Let $f$ be a continuous homogeneous functional of a degree $\gamma$ defined
on convex sets. From Theorem \ref{th2}, we get
$$
k^{\frac{\gamma}{\alpha}}f(C_k(\pi_{\alpha,\nu})) \;\to f(C_{\infty})\;\;\;a.s.
$$

In particular, if $f(A)$ is the surface Lebesgue measure, i.e.
$f(A)=\lambda^{d-1}(\partial A)$, then \begin{equation}
 f(C_k(\pi_{\alpha,\nu})) \sim \frac{f(C_{\infty})}{k^{\frac{d-1}{\alpha}}}.
\end{equation}

But if $f(A) = \lambda^{d}(A),$ then \begin{equation} f(C_k(\pi_{\alpha,\nu})) \sim
\frac{f(C_{\infty})}{k^{\frac{d}{\alpha}}}. \end{equation} \end{remark}

\section{Proof}
\subsection{Auxiliary lemmas}

Let $\eta_1,\eta_2,\dots $ be i.i.d. random variables with the
standard exponential distribution, so that $a={\rm
E}\eta_1={\rm Var}\eta_1=1.$ Define the sums
\begin{equation}
\Gamma_n=\eta_1+\cdots+\eta_n.
\label{eq:defgamman}
\end{equation}
By the law of the iterated logarithm there exists an a.s. finite
random variable $\kappa$  with values in $\bbn$ such that for $ n\geq
\kappa$ 
\begin{equation}\label{LIL}
|n^{-1}\Gamma_n-1|<2\sqrt{\frac{\lln n}{n}} \quad a.s..\end{equation} 

Consider a function $h(z)=z^{-1/\alpha}.$ 
If $|z-1|\le 1/2$ then
$$|h(z)-h(1)|\le L_{\alpha}|z-1|,\ L_{\alpha}<\infty.$$ \ 
Let $n'=\min\left\{n|\
2\sqrt{\frac{\lln n}{n}}<1/2\right\}.$ If $n\ge \max(n' ,\kappa)$ then
$$|h(n^{-1}\Gamma_n)-h(1)|\le 2L_{\alpha}\sqrt{\frac{\lln n}{n}}$$
and, therefore, for $n\geq \kappa = \kappa(\omega)$ 
\begin{equation}\label{LIL1} \left|\Gamma_n^{-1/\alpha} -
n^{-1/\alpha}\right| \le 2L_{\alpha}\frac{\sqrt{\lln{n}}}{n^{1/\alpha
+1/2}} \quad a.s..\end{equation}

We call the {\em configuration} any countable set of points from
$\bbr^d$ such that for any $\delta>0$ there are a finite number of
points belonging to the set that lie outside the ball $\{x\mid\
|x|\le\delta \}.$ So the point ${\mathbf 0}$ is the limit point of
any configuration. We call a configuration $\varsigma$ {\em
non-unilateral} if all the convex hulls, $C_k=C_k(\varsigma),\
k=1,2,\dots,$ generated by $\varsigma$ contain ${\mathbf 0}$ as an
interior point. It is evident that under the conditions of Theorem
\ref{SHAPE} the random measure $\pi_{\alpha,\nu}$ is almost surely supported by a
non-unilateral configuration $\varsigma.$

Denote by ${\rm int}(\varsigma)$ the set of the interior
  points of $\varsigma$, i.e.
 $${\rm int}(\varsigma)=\{x\; |\;  x\in \varsigma,\ x \notin \partial
 C_1(\varsigma)\}.$$

\begin{lemma} \label{LEMMA1vil}
Let $\varsigma_1,\ \varsigma_2\in \mathcal{K}$ be such that
$\varsigma_1\subset \varsigma_2$, then for all 
$ k\in \bbn$, 
\begin{equation}\label{INSERT} C_{k}(\varsigma_1) \subset C_{k}(\varsigma_2).\end{equation}
 \end{lemma}
 \begin{proof}
   It is trivial that $C_{1}(\varsigma_1) \subset C_{1}(\varsigma_2)$. 

Note that if $x$ is an interior point of  $C_{1}(\varsigma_1)$, i.e. 
$x\in {\rm int}( C_{1}(\varsigma_1))$,  then  $x$ is also an interior point of $C_1(\varsigma_2)$, 
therefore 
$$  {\rm int}(C_{1}(\varsigma_1)) \subset {\rm int}(C_{1}(\varsigma_2)) $$
and this implies that 
$$ C_2(\varsigma_1) =  C_{1}( {\rm int}(C_{1}(\varsigma_1)) ) \subset C_1( {\rm int}(C_{1}(\varsigma_2)) ) 
= C_2(\varsigma_2)$$
By induction, the lemma is proved .  
 \end{proof}

\begin{lemma} \label{LEMMA1} Let  $\mathcal{ K}$ be the set of non-unilateral
configurations such that no $d+1$ points lie on the same
hyperplane. Let $\varsigma,\ \varsigma'\in\mathcal{ K}$ be such that $\varsigma' \subset \varsigma$ and  $\#(\varsigma'
\setminus \varsigma)=m < \infty.$ Then we have, for all  $ k\in \bbn$,
\begin{equation}\label{ASSERTION} C_{k+m}(\varsigma) \subset C_{k}(\varsigma') \subset
C_{k}(\varsigma)\end{equation} \end{lemma}

\begin{proof}
  
  Since $\varsigma,\ \varsigma'\in\mathcal{K}$ and ${\mathbf 0}$ is the only limit
  point of both configurations all $C_k(\varsigma),\ k=1,2,\dots,$ are
  polyhedrons.  
Note that for all $k,l\ge 1$ \begin{equation}\label{REC1}
C_{k+1}(\varsigma)=C_1({\rm int}(C_k(\varsigma)\cap\varsigma)) \end{equation} and
\begin{equation}\label{REC2} C_{k+l}(\varsigma)=C_k({\rm int}(C_l(\varsigma)\cap\varsigma)).
\end{equation}

First, let $m=1.$ Note that the inclusion  $C_1(\varsigma') \subset
C_1(\varsigma)$ follows directly from the relation  $\varsigma' \subset \varsigma.$
Denote $ \{a\} =\varsigma
 \setminus \varsigma'.$ Consider two possible cases $a \notin C_1(\varsigma')$
  and  $a \in C_1(\varsigma')$ one after another.

Let $a \notin C_1(\varsigma').$ In this case $a\in\partial C_1(\varsigma)$
i.e. $C_1(\varsigma')\neq C_1(\varsigma).$ It implies that ${\rm
int}(\varsigma)\subset\varsigma'.$ Utilising (\ref{REC1}) under $k=1$ yields
$C_2(\varsigma') \subset C_1(\varsigma).$ Since the inclusion $C_1(\varsigma')
\subset C_1(\varsigma)$ is obvious we conclude that (\ref{ASSERTION})
holds for $k=m=1.$

Further, let us make use of the induction by $k.$ Assume that
(\ref{ASSERTION}) holds for $m=1$ and all $k\leq n$ and show that
then it holds  for $m=1$ and $k=n+1.$ By the induction assumption
we have

  \begin{equation} \label{INDUCTION}
  C_{n+1}(\varsigma) \subset C_{n}(\varsigma') \subset C_{n}(\varsigma).
\end{equation} Since
$${\rm int}(C_n(\varsigma')\cap\varsigma')={\rm int}(C_n(\varsigma')\cap\varsigma)$$
we obtain, taking into account (\ref{REC1}),
$$C_{n+l}(\varsigma')=C_1({\rm int}(C_n(\varsigma')\cap\varsigma)).$$
From the right hand side inclusion of (\ref{INDUCTION}), it follows
that   $ C_{n+1}(\varsigma') \subset C_{n+1}(\varsigma).$  Further, from the
left hand side inclusion of (\ref{INDUCTION}) we conclude that
$${\rm int}(C_{n+1}(\varsigma)\cap\varsigma)={\rm int}(C_n(\varsigma')\cap\varsigma').$$
Applying (\ref{REC2}) yields $C_{n+2}(\varsigma) \subset C_{n+1}(\varsigma').$
Thus, (\ref{ASSERTION}) holds for $k=n+1$ and  $m=1,$ i.e. the case $a
\notin C_1(\varsigma')$ is exhausted.

If  $a \in C_1(\varsigma'),$ then there exists an integer $n_0$ such
that 
\begin{align*}
   C_{n}(\varsigma') & =C_{n}(\varsigma'),\quad n=1,2,\dots,n_0 \\
  C_{n_0+1}(\varsigma') & \not =C_{n_0+1}(\varsigma').
\end{align*}
 Furthermore, $a \notin
C_{n_0+1}(\varsigma').$ Obviously, the relations (\ref{ASSERTION}) are
trivial for $m=1$ and $n=1,2,\dots,n_0.$ Hence, it remains to apply the
above argument to the configurations $C_{n_0+1}(\varsigma)\cap\varsigma$ and
$C_{n_0+1}(\varsigma')\cap\varsigma'.$ Thus, the lemma is proved for all $k$
and $m=1.$

Now, let $m>1,$ i.e. $\varsigma\setminus\varsigma'=\{a_1,\dots,a_m\}.$

Consider the configurations
\begin{align*}
 & \varsigma_0  =\varsigma,  \qquad  
\varsigma_1 = \varsigma \setminus \{a_1\}, \\ 
&\varsigma_2 = \varsigma \setminus \{a_1,a_2\},   \ldots,\\ 
& \varsigma_m =\varsigma\setminus\{a_1,\dots,a_m\}=\varsigma'.
\end{align*}
 
Note that the neighbouring configurations differ by a single point. So, one may
apply (\ref{INDUCTION}). Applying it yields
$$ C_{k+m}(\varsigma) \subset C_{k+m-1}(\varsigma_1) \subset
C_{k+m-2}(\varsigma_2) \subset\ldots \subset C_k(\varsigma_m)=C_k(\varsigma')
\subset C_k(\varsigma). $$ The lemma is proved.
\end{proof}

\subsection{Proof of Theorem \ref{SHAPE}}

\begin{lemma}
\label{LEMMA3} Let $e^{(i)}$ for $i=1,\ldots,l$ with $l\geq d+1,$ be unit
vectors such that $${\rm cone}\ \{e^{(1)},\ldots,e^{(l)}\} =
\bbr^d.$$  If
$A=\{\nu_1^{1/\alpha}e^{(1)},\dots,\nu_l^{1/\alpha}e^{(l)}\}$
 is extreme then there exists  $r>0$ and  $\varepsilon$ depending only on $A$ and on the 
dimension $d$ such that the set $\{r_1\nu_1^{1/\alpha}e^{(1)},\dots,r_l\nu_l^{1/\alpha}e^{(l)}\}$ is extreme 
for all $(r_1,\ldots,r_n)$  such  that   $|r_i/r-1|<\varepsilon,\
i=1,\dots,l$. 
 \end{lemma}
\begin{proof}
  The proof of this lemma  is evident. 
\end{proof}

\begin{proof}[Proof of Theorem \ref{SHAPE}]
  Let us label the points lying on the ray $\mathcal{ L}_i$ in the
  descending order of their norms. So, we have the sequence $
  x_1^{(i)}, x_2^{(i)},\ldots\,$ such that a.s. $
  |x_1^{(i)}|>|x_2^{(i)}|>\cdots.$ Obviously, the sequences
  $\{x_n^{(i)}, n\in \bbn\},\ 1\le i\le l,$ are jointly independent.
  Furthermore, from (\ref{CNTRLMSR})
$$\{|x_n^{(i)}|\} \eqd \{\nu_i^{1/\alpha}\Gamma_n^{-1/\alpha}\}$$
where $\Gamma_n$ is defined as in (\ref{eq:defgamman}).

Let $\epsilon >0$.  According to (\ref{LIL1}), there exists
$n_0=n_0(\omega)$ such that for all $i=1,\ldots,l$ and all $n\ge
n_0$ \begin{equation} \label{INEQUALITY1} \left|x_n^{(i)}
    -\nu_i^{1/\alpha} n^{-1/\alpha}\right| \leq 2L_{\alpha}n^{-1/\alpha
    -1/2}\sqrt{\lln{n}}. \end{equation} and
$$2L_{\alpha}n^{-1/2}\sqrt{\lln{n}}<\varepsilon.$$

Let the configuration $\varsigma'$ is formed by the points $x_n^{(i)},\ n\ge
n_0,\ i=1,\dots,l$ i.e.
$$\varsigma'=
\bigcup_{i=1}^l \{x_{n_0}^{(i)}, x_{n_0+1}^{(i)},\ldots\}.
$$

Consider for all $k\leq 1$
$$A_k^+=\{(n_0+k-1)^{-1/\alpha}(1+\varepsilon_k) \nu_1^{1/\alpha} e^{(1)},\ldots,(n_0+k-1)^{-1/\alpha}(1+\varepsilon_k)\nu_l^{1/\alpha}e^{(l)}\}$$
and
$$A_k^-=\{(n_0+k-1)^{-1/\alpha}(1-\varepsilon_k) \nu_1^{1/\alpha}e^{(1)},\ldots,(n_0+k-1)^{-1/\alpha}(1-\varepsilon_k)\nu_l^{1/\alpha}e^{(l)}\}$$
where
$$\varepsilon_k=2L_{\alpha} \sqrt\frac{\lln (k+n_0-1)}{k+n_0-1}.$$
By virtue of (\ref{INEQUALITY1}) the points
$x_{n_0}^{(1)},\ldots,x_{n_0}^{(l)}$ hit the layer $C(A_1^+)\setminus
C(A_1^-).$ Then by Lemma \ref{LEMMA3} the convex hull $C_1(\varsigma')$ is
the polyhedron and
$${\rm ext}\ C_1(\varsigma')=\{x_{n_0}^{(1)},\ldots,x_{n_0}^{(l)}\}.$$
Similarly, the set
$$\{x_{n_0+1}^{(1)},\ldots,x_{n_0+1}^{(l)}\}\subset
C(A_2^+)\setminus C(A_2^-)$$ and, therefore, it is extreme, i.e.
$${\rm ext}\ C(\{x_{n_0+1}^{(1)},\ldots,x_{n_0+1}^{(l)}\})=
\{x_{n_0+1}^{(1)},\ldots,x_{n_0+1}^{(l)}\}.$$ It is evident that
$${\rm ext}\ C_2(\varsigma')=\{x_{n_0+1}^{(1)},\ldots,x_{n_0+1}^{(l)}\}.$$
Continuing in this way we obtain at the $k$-th convex hull $C_k(\varsigma')$
such that $${\rm ext}\ C_k(\varsigma')=
\{x_{n_0+k-1}^{(1)},\ldots,x_{n_0+k-1}^{(l)}\}\subset
C(A_k^+)\setminus C(A_k^-).$$ The last inclusion implies that
$$ d_H\left((k+n_0-1)^{1/\alpha}C_k(\varsigma'),\,C_{\infty}\right)
\le\varepsilon_k, $$ where, we remind, $C_{\infty}$ is the convex hull
generated by
$A=\{\nu_1^{1/\alpha}e^{(1)},\ldots,\nu_l^{1/\alpha}e^{(l)}\}.$ From
(\ref{ASSERTION}) it follows that
$$C_{k+m}(\varsigma')\subset C_{k+m}(\pi_{\alpha,\nu})\subset C_k(\varsigma'),\quad \textrm{with } m=(n_0-1)l.$$ Therefore, $${\rm ext}\ C_{k+m}(\pi_{\alpha,\nu})\subset
C(A_k^+)\setminus C(A_{k+m'}^-),\quad \textrm{with }
m'=(n_0-1)(l-1).$$ So, for all sufficiently large $k$

$$ d_H\left((k+m)^{1/\alpha}C_{k+m}(\pi_{\alpha,\nu}),\,C_{\infty}\right) \le
2\varepsilon_k.$$ Since $m$ is fixed the theorem follows.
\end{proof}
\subsection{Proof of Theorem \ref{th2}}
  Let $\epsilon$ be an arbitrary positive real.  Hereafter, we denote
  $ A^{(\epsilon)}$ the $\epsilon$-neighbourhood of a set $A$,
 $$ A^{(\epsilon)}=\{ x  \; : \;  d(x,A)<\epsilon\}.$$

 \vskip 0.5cm Let $A_1 = A\cap \partial C(A) \stackrel{def}{=} \{
 \nu_j^{1/\alpha}e^{(j)} , j\in J\}$, the set $A_1$ is extreme.  From the
 process $\pi_{\alpha,\nu}$ , we construct a new p.p.p.  $\pi_1 $
 obtained by deleting all the points on the rays $\mathcal{ L}_j=\{x|\
 x=te^{(i)},\ t>0\}$, $j\in J$.  By Lemma \ref{LEMMA1vil}, we have for
 all $n\in \mathbb{N}$
 \begin{equation}
   \label{eq:d2}
   C_n(\pi_1) \subset  C_n(\pi_{\alpha,\nu})
 \end{equation}Moreover,   $A_1$ is extreme and $C_{\infty}=C(A) =C(A_1)$, 
 thus Theorem \ref{SHAPE} ensures 
 the  convergence 
 \begin{equation}
   \label{eq:d3}
   d_H( n^{1/\alpha} C_n(\pi_1) ,  C_{\infty}) \to  0 \; {\rm a.s.}    
 \end{equation}
 From (\ref{eq:d2}) and (\ref{eq:d3}), it exists $n_1\in \mathbb{N} $
 such that for all $n>n_1$
 \begin{equation}
   \label{eq:d6}
   C_\infty\subset n^{1/\alpha}C_n(\pi_1)^{(\epsilon)}  \subset n^{1/\alpha}C_n(\pi_{\alpha,\nu})^{(\epsilon)}.   
 \end{equation}
 It is easy to see that there exists $\tilde \nu_i, \; i\in I$ such
 that the set.  
$$ A_2 = \{    \nu_j^{1/\alpha}e^{(j)} ,\  j\in J\ ; \; ; \tilde\nu_i^{1/\alpha}e^{(i)} ,\  i\in  \{1,\ldots,l\} \setminus J \} $$
is extreme and satisfies the following relation
\begin{equation}
  \label{eq:d3p}
  C_\infty \subset C(A_2)  \subset C_\infty^{(\epsilon)}.   
\end{equation}

From $\pi_{\alpha,\nu}$ we construct a second p.p.p.  $\pi_2 $ by
adding the independent point processes $(\tilde\pi_i)_{i\in J}$
verifying the following conditions
\begin{itemize}
\item the $(\tilde\pi_i)_{i\in J}$ are independent of
  $\pi_{\alpha,\nu}$ ;
\item for each $i\in J$ the spectral measure of $\tilde\pi_i$ is
  supported by $\mathcal{ L}_i$ and the intensity measure is
  $\tilde\mu_i(A) = (\tilde \nu_i -\nu_i) \alpha \int_A r^{-\alpha-1}
  {\rm d }\ r$.
\end{itemize}

According to Theorem \ref{SHAPE} we have
\begin{equation}
  \label{eq:d4}
  d_H( n^{1/\alpha} C_n(\pi_2) ,  C(A_2)  \to 0 \; {\rm a.s.}
\end{equation}
and using Lemma \ref{LEMMA1vil}, we get, for all $n\in\bbn$,
\begin{equation}
  \label{eq:d5}
  C_n(\pi_{\alpha,\nu})  \subset C_n(\pi_2).
\end{equation}

From (\ref{eq:d3p}),(\ref{eq:d4}) and (\ref{eq:d5}), it exists $n_2\in
\mathbb{N} $ such that for all $n>n_2$
\begin{equation}
  \label{eq:d7}
  n^{1/\alpha} C_n(\pi_{\alpha,\nu} )  \subset n^{1/\alpha}C_n(\pi_2) \subset C(A_2)^{(\epsilon)}  \subset C_\infty^{(2\epsilon)}  .
\end{equation}

According to (\ref{eq:d6}) and (\ref{eq:d7}), for all $n\geq
\max(n_1,n_2)$, we have
$$ n^{1/\alpha} C_n(\pi_{\alpha,\nu}) \subset  C_\infty^{(2\epsilon)} \qquad {\rm and} \qquad   C_\infty  \subset  n^{1/\alpha} C_n(\pi_{\alpha,\nu})^{(2\epsilon)}.$$
By definition of $d_H$, this means $$d_H(n^{1/\alpha}
C_n(\pi_{\alpha,\nu}) , C_\infty) \le 2\epsilon,$$ and we get
(\ref{CV22}).
\section{Simulation and conjectures}
We investigate using some simulations the limit shape and the asymptotic behaviour of basic functionals in the case of  continuous spectral measure.  Hereafter we consider the example of the uniform distribution  as  spectral measure. 

The point processes  $\{ x^{(j)} \;, \; j\in \mathbb{N}\}$  are simulated using the representation (\ref{ppsim}). 
Let $C_{1,n}$ be  the convex hull generated by the first $n$ points
$x^{(1)},\ x^{(2)},\dots,x^{(n)} $ and  $\kappa_{n,1}=\displaystyle\min_{x\in\partial C_{1,n}}\ |x|.$
Since the points of the simulated  p.p.p.
are ordered by their distances from the origin, it is evident that 
 $$C_{1,n'}=C_1 
\text{ with }  n'=\min\{n  \; : \; \kappa_{n,1}>|x^{(n+1)}|\}. $$ 
This fact is  used to construct the successive convex hulls  $(C_{k})_{k\in \mathbb{N}}$.

\begin{figure}[htbp]
\centering%
\includegraphics[scale=.5]{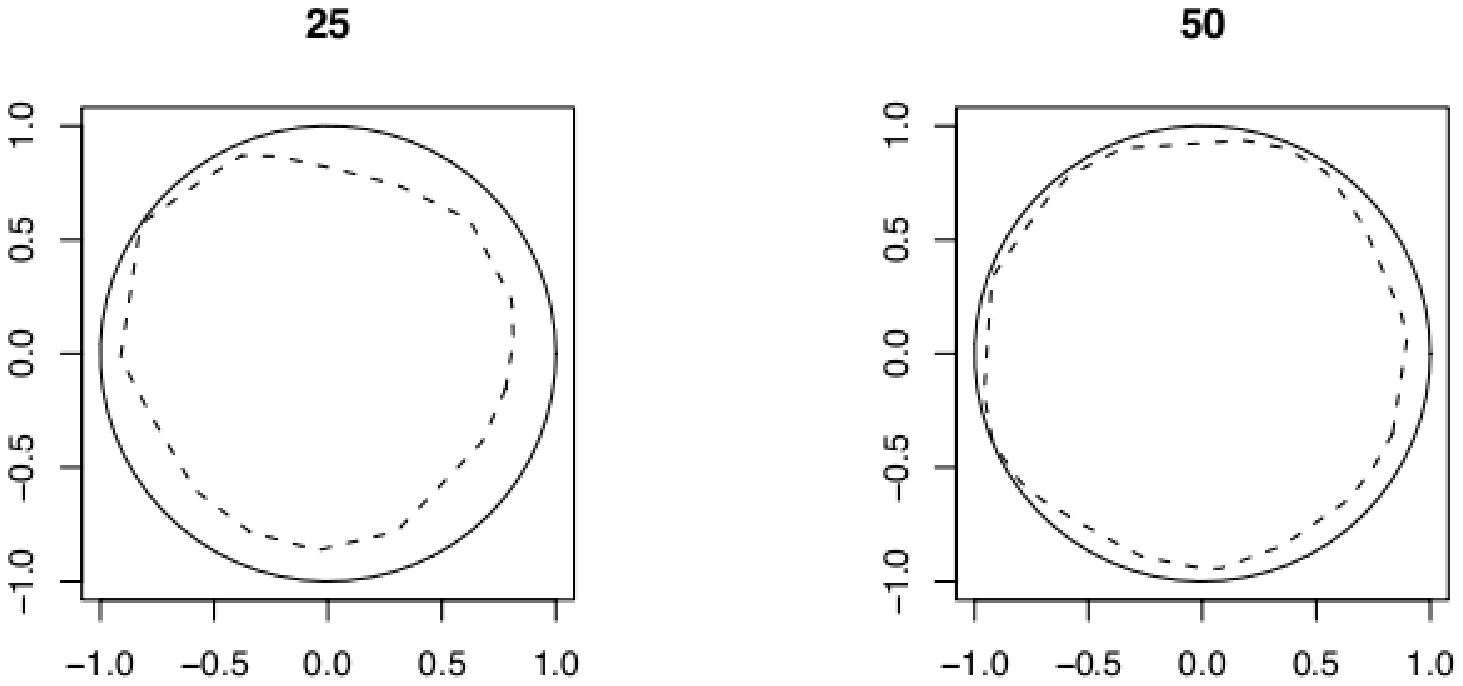}\hskip .3cm 
\includegraphics[scale=.5]{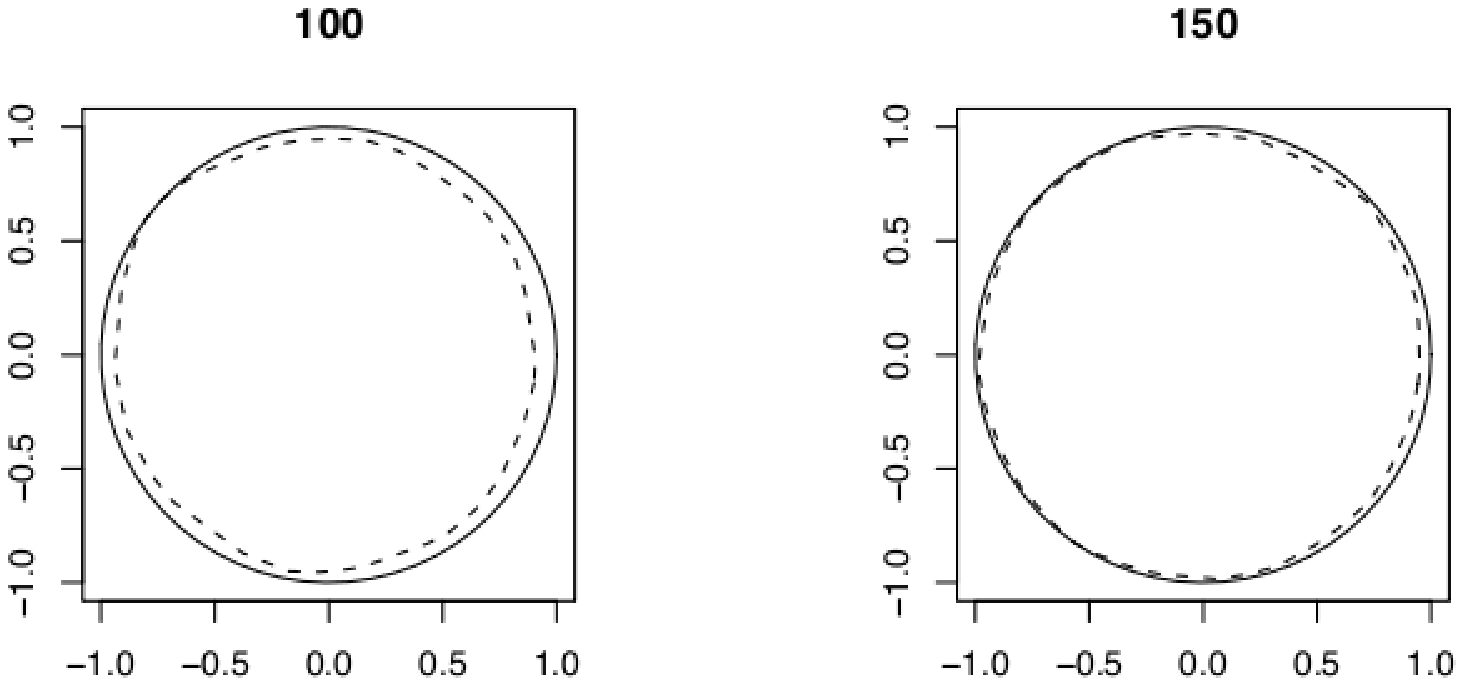}
  \caption{[dotted line] the normalized simulated shapes of $\hat C_{25}$, $\hat C_{50}$, $\hat C_{100}$
  and $\hat C_{150}$ in the case  $\alpha=3/2$ and the spectral measure
  is uniform. [solid line]  the unit circle. }
 \label{fig:uniform2}\end{figure}

 Figures \ref{fig:uniform2} gives an impression about
the behaviour of the peels.  The observed closeness of the peels to the unit circle also support our conjecture about the existence of the limit shape that is expected to be a circle.

\begin{figure}[htbp]
  \centering
  \hbox{
   \includegraphics[height=4cm,width=6cm,angle=270]{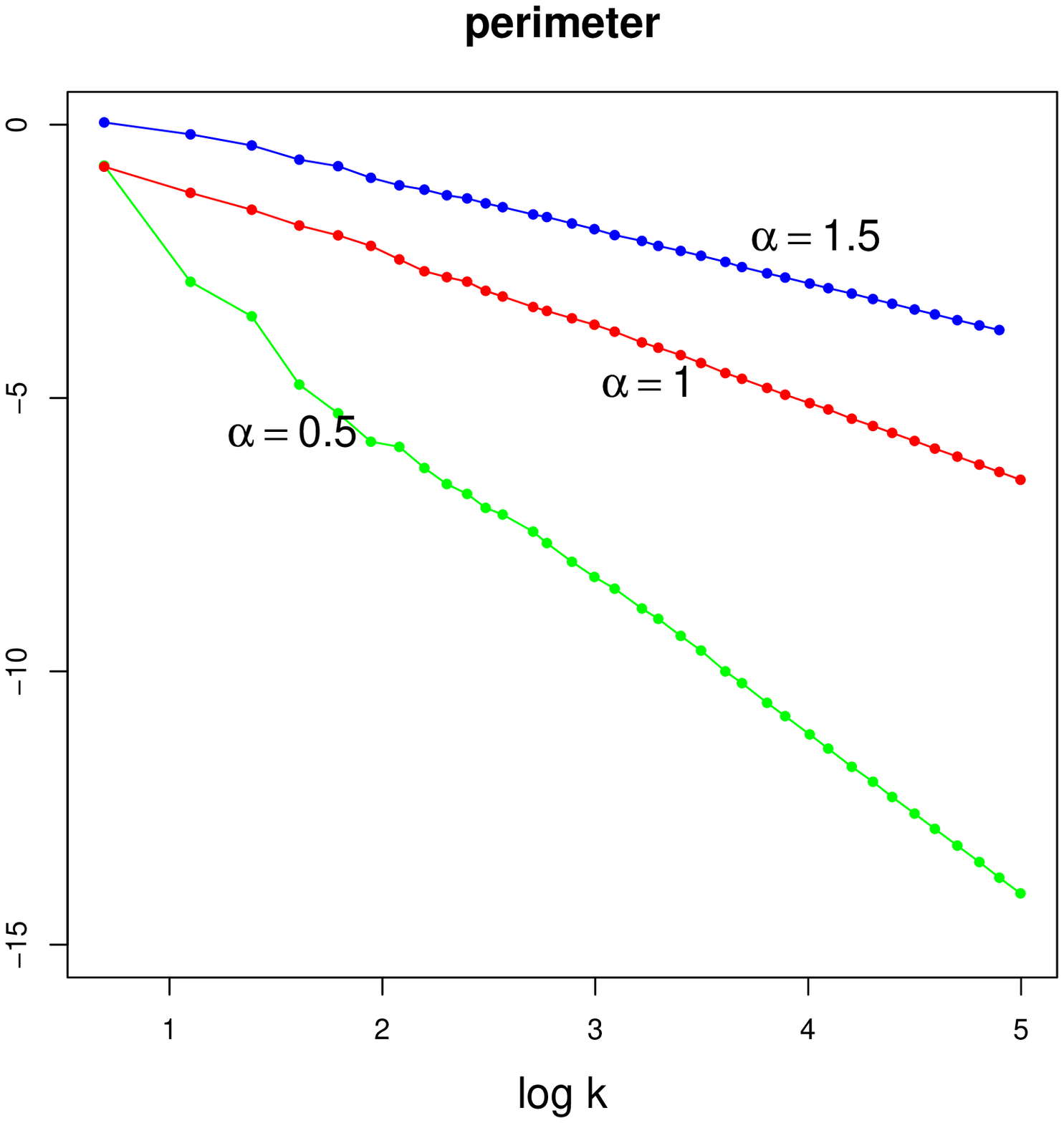}
  \hskip -.8cm \includegraphics[height=4cm,width=6cm,angle=270]{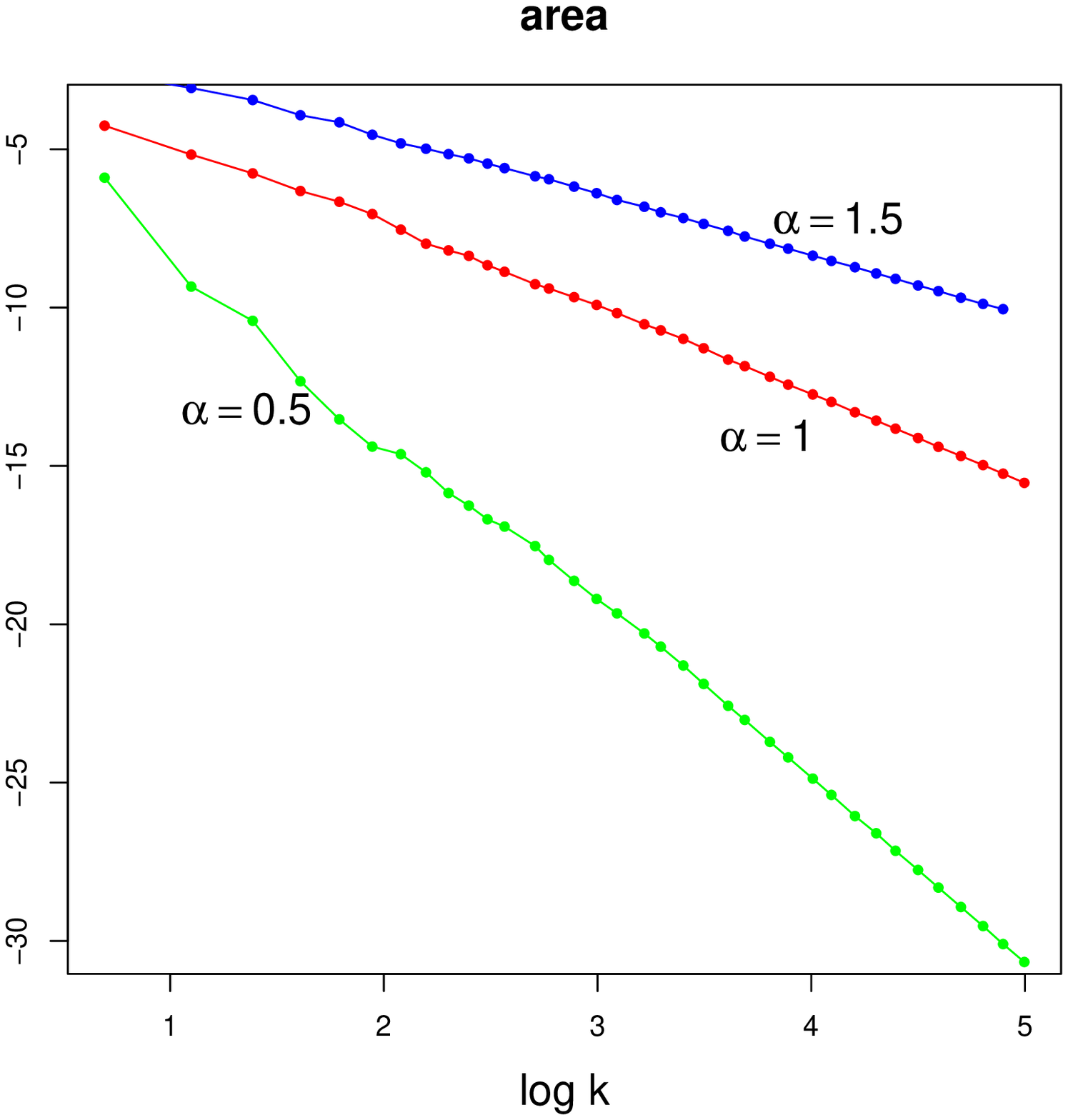}
  \hskip -.8cm   \includegraphics[height=4cm,width=6cm,angle=270]{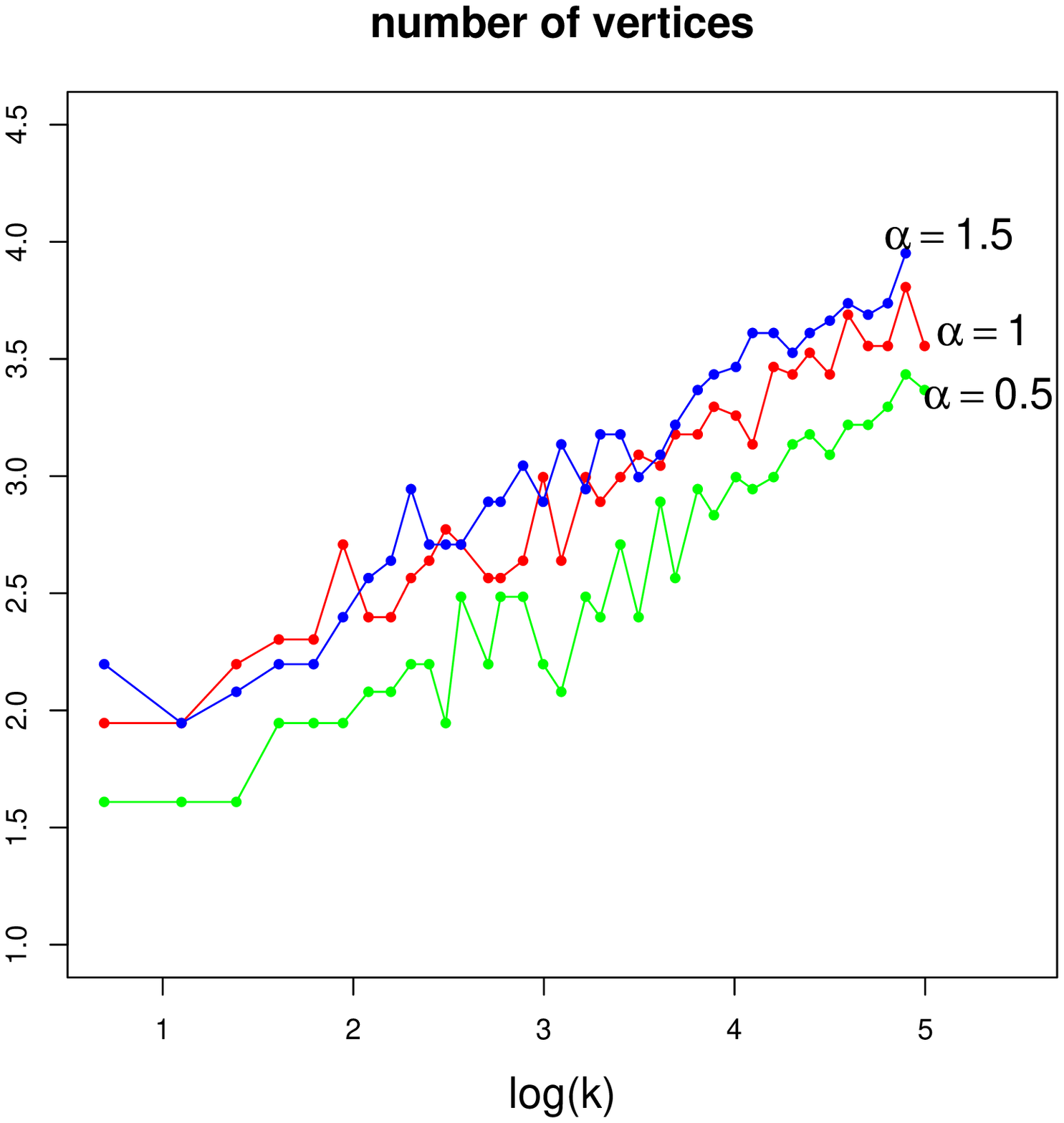}
 }
 \vskip -.5cm
  \includegraphics[height=6cm,width=8cm]{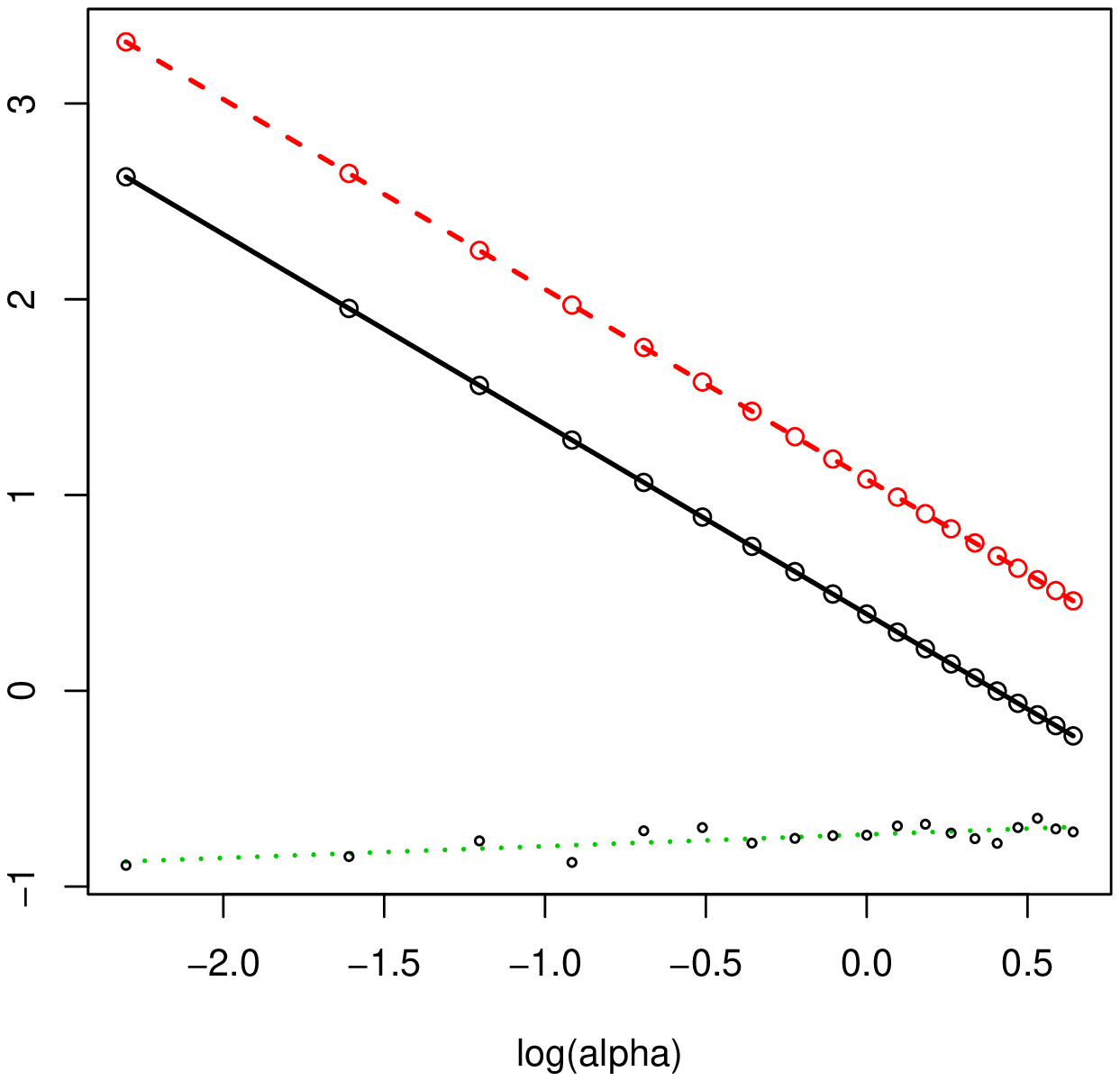}
 \vskip -.5cm
 \caption{[TOP] Log-log representation of the  values of   $\mathcal{ L}(C_{k}),\ \mathcal{ A}(C_{k})$ and $\mathcal{ N}(C_{k})$  as function of $k$. The functionals are calculated on simulated p.p.p.  for different values of $\alpha=0.5;1.0;1.5 $ and the uniform distribution as spectral measure. 
 [BOTTOM] The estimated values of the logarithm of exponents defined in (\ref{hypo}) versus $\ln(\alpha)$  and the best linear fittings  $\ln\hat{\gamma}_{a}=-0.97\ln\alpha+\ln 2.95$,  
   $\ln\hat {\gamma}_{l}=-0.97\ln\alpha+\ln 1.48$ and $ \ln\hat{\gamma}_{n}=0.06\ln\alpha+\ln 0.48$.  Three exponents are estimated for each $\alpha$ on $1000$ independent replications.
 }
     \label{fig:rate}
 \end{figure}

It is of great interest to get impression about a possible
behaviour of such basic functionals of the convex polygons $C_k,\
k=1,2,\dots,$ as the perimeter $\mathcal{ L},$ the area $\mathcal{ A}$ and
the total number of vertices $\mathcal{N}.$ 
It seems evident that   $ \mathcal{ L}(C_k)$ and $\mathcal{ A}(C_k)$ tend to zero as $k\to \infty$. Intuitively, we expect that  $\mathcal{
N}(C_k)\to\infty$ as $k\to\infty$.  
 Figure \ref{fig:rate}[top]  represents  the logarithm of those functionals as function of $\log(k)$, calculated on simulated p.p.p. for different values of $\alpha$. 
The observed closeness of the points to straight lines makes it reasonable to expect that in a sense 
\begin{equation} \label{hypo} \mathcal{ L}(C_k) \asymp
k^{-\gamma_l},\  \mathcal{ A}(C_k) \asymp k^{-\gamma_a},\  \mathcal{ N}(C_k)\asymp k^{\gamma_n}\end{equation} with  $\gamma_{l}$ and $\gamma_{a}$ depending on $\alpha$ whereas it seems that $\gamma_{n}$ does not dependent on $\alpha$. 
 
 The next step consists in estimating  the
exponents and possibly the dependence  on $\alpha.$  Using independent replications of p.p.p., we estimate the three exponents defined in (\ref{hypo}) for different values  of $\alpha$.  Figure \ref{fig:rate} [bottom] represents the logarithm of the estimated exponents versus  $\log(\alpha)$.  For the three cases, the linear approximation seems reasonable. According to the estimated coefficients of the straight lines  (see the equations in the caption of  Figure \ref{fig:rate}), it looks very credible that   the true values are 
 \begin{equation}
\gamma_{l}=\frac{3}{2\alpha},\ \gamma_{a}=\frac{3}{\alpha}\text{ and } \gamma_{n}=\frac{1}{2}.
 \label{est.slo}
\end{equation}

After the k-th iterative step of the peeling procedure, the number of deleted points  should be  the order of 
 $\sum_{j=1}^{k} \mathcal{ N}(C_j) \asymp  k^{\frac 32}$ and   
 $$\rho_{k}=\max_{x \in C_{k}} |x| \asymp \left(  k^{\frac 32} \right)^{-1/\alpha}   $$
Moreover we can expect that $
 d_H\left(\rho_{k}^{-1} C_k(\pi_{\alpha,\nu}),\; C_{\infty}\right)$ converges to zero.   

Using the arguments of  Remark  \ref{az}, this convergence  would lead to  $ \mathcal{ L}(C_k) \asymp \rho_{k} \asymp k^{-3/(2\alpha)}$ and $ \mathcal{ A}(C_k) \asymp \rho_k^{2}\asymp k^{-3/\alpha}$. 
These convergence rates are in agreement with the estimated values  of $\alpha_{l}$ a,d $\alpha_{a}$  obtained in (\ref{est.slo}).

\end{document}